# Existence and convergence theorems for monotone generalized $\alpha$–nonexpansive mappings in uniformly convex partially ordered hyperbolic metric spaces and its application

Chang Il Rim, Jong Gyong Kim, Chol-Hui Yun

# Existence and convergence theorems for monotone generalized $\alpha$–nonexpansive mappings in uniformly convex partially ordered hyperbolic metric spaces and its application

Chang Il Rim, Jong Gyong Kim, Chol-Hui Yun

Faculty of Mathematics, **Kim Il Sung** University, Pyongyang, DPR Korea

**Abstract** In this paper, we generalize the existence result in [14] and prove convergence theorems of the iterative scheme in [12, 16] for monotone generalized $\alpha$–nonexpansive mappings in uniformly convex partially ordered hyperbolic metric spaces. And we also give a numerical example to show that this scheme converges faster than the scheme in [14] and apply the result to the integral equation.

## 1. Introduction and preliminaries

In many literatures, several iterative schemes for fixed points of nonexpansive mappings have been studied.

In 1953, Mann([3]) proposed the following iterative scheme for nonexpansive mappings.

$$\begin{cases} x = x_1 \in C \\ x_{n+1} = (1-a_n)x_n + a_n T x_n \end{cases}, n \in \mathbf{N} \qquad (*)$$

, where $\{a_n\}, \{b_n\}, \{c_n\} \subset (0,1)$.

In 2016, Sahu([12]) and Thakur([16]) introduced the following iterative scheme for nonexpansive mappings in uniformly convex Banach spaces and proved that this scheme converges to a fixed point of a contraction mapping faster than all the known iterative schemes.

$$\begin{cases} x = x_1 \in C \\ x_{n+1} = (1-a_n)Tz_n + a_n Ty_n \\ y_n = (1-b_n)z_n + b_n Tz_n \\ z_n = (1-c_n)x_n + c_n Tx_n \end{cases}, n \in \mathbf{N} \quad (**)$$

, where $\{a_n\}, \{b_n\}, \{c_n\} \subset (0, 1)$.

On the other hand, in recent years some new generalized nonexpansive mappings have been introduced and fixed point theorems and convergence theorems for such mappings have been studied.

In 2008, Suzuki([15]) defined a new kind of generalized nonexpansive mappings, ($C$)-mapping.

**Definition 1** ([15]) Let $T$ be a mapping on the subset $C$ of a Banach space. $T$ is said to satisfy the condition ($C$)(or $T$ is a ($C$)-mapping) if it satisfies the following property for all $x, y \in C$.

$$\frac{1}{2}\|x - Tx\| \leq \|x - y\| \Rightarrow \|Tx - Ty\| \leq \|x - y\|.$$

Since then, several existence theorems and convergence theorems for a fixed point of ($C$)-mappings have been developed. In 2019, [5] introduced (**) for ($C$)-mappings in uniformly convex Banach spaces to prove the convergence theorems to a fixed point and gave a numerical example to show that this scheme converges faster than other schemes.

Also, fixed point theorems and convergence theorems for monotone generalized nonexpansive mappings in partially ordered Banach spaces have been studied.

In 2017, [14] defined another new kind of monotone generalized nonexpansive mapping, monotone generalized $\alpha$-nonexpansive mappings, and obtained convergence and existence theorems of the scheme (*) in uniformly convex partially ordered Banach spaces.

All results mentioned above are all obtained in Banach space. However, reality problems are not all given in Banach spaces. For this reason, scientists defined a more general space than Banach space, hyperbolic metric space, and studied fixed point theorems and convergence theorems for nonexpansive mappings and its generalization on such space.

In this paper, we prove the fixed point theorem and convergence theorems of the iterative scheme (**) for monotone generalized –nonexpansive mappings and give a numerical example to show that (**) converges faster than (*). Also we give its application to an integral equation.

Here are some basic concepts.

**Definition 2** ([9]) Let $(X, \rho)$ be a metric space. Assume that the mapping $H: X \times X \times [0,1] \to X$ satisfies the following conditions for all $u, v, w, z \in X$ and $\beta, \gamma \in [0,1]$.

1. $\rho(z, H(u, v, \beta)) \leq (1-\beta)\rho(z, u) + \beta\rho(z, v)$;
2. $\rho(H(u, v, \beta), H(u, v, \gamma)) = |\beta - \gamma|\rho(u, v)$;
3. $H(u, v, \beta) = H(v, u, 1-\beta)$;
4. $\rho(H(u, z, \beta), H(v, w, \beta)) \leq (1-\beta)\rho(u, v) + \beta\rho(z, w)$;

Then $(X, \rho, H)$ is said to be a hyperbolic metric space.

We write $H(u, v, \beta) = (1-\beta)u \oplus \beta v$ as well.

Every linear normed space is a hyperbolic metric space and Hadamard manifold([2]) and Hilbert open unit ball equipped with the hyperbolic metric([4]) are examples of hyperbolic metric spaces which are not normed spaces.

Let $Y$ be a subset of hyperbolic metric space $X$. $Y$ is said to be convex if $(1-\beta)u \oplus \beta v \in Y$ holds, for any $u, v \in Y$ and $\beta \in [0,1]$.

Let $(E, \rho, \leq)$ be a partially ordered hyperbolic metric space equipped with the partial order $\leq$. Then for any $a, b \in E$, we denote order intervals as follows.

$$[a, \to) := \{u \in E : a \leq u\}, \quad (\leftarrow, b] := \{u \in E : u \leq b\}$$

Throughout this paper, we assume that these order intervals are closed and convex.

**Definition 3** ([6]) Let be a hyperbolic metric space. For any $a \in E$ and $r > 0, \varepsilon \in (0, 2]$, we define the mapping $\delta: (0, \infty) \times (0, 2] \to (0, 1]$ as follows.

$$\delta(r, \varepsilon) := \inf\left\{1 - \frac{1}{r}\rho\left(\frac{1}{2}u \oplus \frac{1}{2}v, a\right) : \rho(u, a) \leq r, \rho(v, a) \leq r, \rho(u, v) \geq r\varepsilon\right\}.$$

If $\delta(r, \varepsilon) > 0$ holds for any $r > 0, \varepsilon \in (0, 2]$, then $E$ is said to be uniformly convex. And the mapping $\delta$ is called modulus of uniform convexity. $\delta$ is said to be monotone if it decreases with $r$ for a fixed $\varepsilon$.

**Definition 4** ([14]) Let $\{x_n\}$ be a bounded sequence of a hyperbolic space $E$. For $x \in E$, we set the asymptotic radius of $\{x_n\}$ relative to $x$ as follows.

$$r(x, \{x_n\}) := \limsup_{n \to \infty} \rho(x_n, x).$$

The asymptotic radius of $\{x_n\}$ relative to $C$ is defined by

$$r(C, \{x_n\}) := \inf \{r(x, \{x_n\}) | \ x \in C\}.$$

The asymptotic center of $\{x_n\}$ relative to $C$ is the set

$$A(C, \{x_n\}) := \{x \in C | \ r(x, \{x_n\}) = r(C, \{x_n\}) \}.$$

It is well known that $A(C, \{x_n\})$ consists only one point with respect to any closed convex subset in complete uniformly convex hyperbolic space with monotone modulus of uniform convexity([10]).

**Definition 5** ([1]) Let $C$ be a nonempty subset of a hyperbolic space $(E, \rho)$. If for any $x \in C$, there exists a bounded sequence $\{x_n\} \in E$ such that

$$\tau(x) = \limsup_{n \to \infty} \rho(x_n, x),$$

then the function $\tau : C \to [0, \infty)$ is called a type function.

We know that every bounded sequence generates a unique type function.

In uniformly convex hyperbolic space, a type function on a closed convex subset is continuous and has a unique minimum point([1]).

**Definition 6** ([8]) Let $(E, \rho)$ be a hyperbolic metric space. A bounded sequence $\{x_n\} \subset E$ is said to $\Delta$–converge to $z \in E$, if $z$ is a unique asymptotic center of every subsequence $\{x_{n_j}\}$ of $\{x_n\}$.

We can see that $\{x_n\}$ $\Delta$-converges to a point $z \in E$, if $z$ is unique and the type function generated by every subsequence $\{x_{n_j}\}$ of $\{x_n\}$ attains its infmum at $z$.

**Definition 7** ([13]) A map $T : C \to C$ is said to satisfy condition ($I$), if there is a nondecreasing function $h : [0, \infty) \to [0, \infty)$ with $h(0) = 0$ and $h(r) > 0$ for any $r > 0$, such that

$$d(x, Tx) \geq h(d(x, F(T))),$$

for all $x \in C$, where $d(x, F(T)) = \inf \{d(x, p) : p \in F(T)\}$.

**Definition 8** ([14]) Let $T: C \to C$ be a mapping on a nonempty subset $C$ of a partially ordered metric space $(E, \rho, \leq)$. $T$ is said to be monotone if for any $x, y \in C$ with $x \leq y$, $Tx \leq Ty$ holds.

$T$ is said to be a monotone generalized $\alpha$-nonexpansive mapping if it is monotone and there exists $\alpha \in [0, 1)$ such that the following property is satisfied for any with $x, y \in C$ with $x \leq y$.

$$\frac{1}{2}\rho(x, Tx) \leq \rho(x, y) \Rightarrow \rho(Tx, Ty) \leq \alpha\rho(Tx, y) + \alpha\rho(x, Ty) + (1 - 2\alpha)\rho(x, y).$$

Every monotone ($C$)-mapping is a monotone generalized $\alpha$-noexpansive mapping. And every monotone generalized $\alpha$-noexpansive mapping is monotone quasi-nonexpansive([14]).

## 2. Existence theorem

In this section, we establish the existence theorem for monotone generalized $\alpha$-noexpansive mappings in uniformly convex partially ordered hyperbolic space.

In a hyperbolic space, (**) is defined as follows.

$$\begin{cases} x = x_1 \in C \\ x_{n+1} = (1 - a_n)Tz_n \oplus a_n Ty_n \\ y_n = (1 - b_n)z_n \oplus b_n Tz_n \\ z_n = (1 - c_n)x_n \oplus c_n Tx_n \end{cases}, n \in \mathbf{N}$$

, where $\{a_n\}, \{b_n\}, \{c_n\} \subset (0, 1)$.

First we give some useful lemmas.

The following lemma shows the property of monotone generalized $\alpha$-noexpansive mappings, which can be proved similarly to lemma 3.8 in [14].

**Lemma 1** Let $C$ be a nonempty subset of a partially ordered hyperbolic space $(E, \rho, \leq)$ and $T: C \to C$ be a monotone generalized $\alpha$-noexpansive mapping on $C$. Then for any $x, y \in C$ such that $x$ and $y$ are comparable, the following inequality holds.

$$\rho(x, Ty) \leq \frac{3 + \alpha}{1 - \alpha}\rho(x, Tx) + \rho(x, y).$$

**Theorem 1** Let $C$ be a nonempty closed convex subset of a uniformly convex partially ordered hyperbolic space $(E, \rho, \leq)$ and $T: C \to C$ be a monotone generalized $\alpha$-noexpansive mapping on $C$. Assume that there exists a bounded sequence $\{x_n\} \subset C$ with $x_n \leq Tx_n (Tx_n \leq x_n)$ for all $n \in \mathbf{N}$, and there exists $x \in C$ such that $x_n \leq x (x \leq x_n)$ for all $n \in \mathbf{N}$. If $\liminf_{n \to \infty} \rho(x_n, Tx_n) = 0$, then $T$ has a fixed point.

**Proof** Assume that $x_n \leq Tx_n$ for all $n \in \mathbf{N}$. Since $\liminf_{n \to \infty} \rho(x_n, Tx_n) = 0$, there exists a subsequence $\{x_{n_k}\}$ of $\{x_n\}$ such that $\lim_{k \to \infty} \rho(x_{n_k}, Tx_{n_k}) = 0$.

$C_j := \{x \in C : x_{n_j} \leq x\}$ is closed convex and from the assumption $C_\infty := \bigcap C_j$ is a nonempty and closed convex subset of $C$. Since $T$ is monotone, for any $y \in C_\infty$, we have $x_{n_j} \leq Tx_{n_j} \leq Ty$ for all $j \in \mathbf{N}$. Hence we have $T(C_\infty) \subset C_\infty$.

Let $\tau: C_\infty \to [0, \infty)$ be a type function generated by $\{x_{n_k}\}$. Then $\tau$ has a unique minimum point $z$. Now we prove that $z$ is a fixed point.

From lemma 1, we have

$$\begin{aligned}\tau(Tz) &= \limsup_{k \to \infty} \rho(x_{n_k}, Tz) \\ &\leq \frac{3+\alpha}{1-\alpha} \limsup_{k \to \infty} \rho(x_{n_k}, Tx_{n_k}) + \limsup_{k \to \infty} \rho(x_{n_k}, z) \\ &= \limsup_{k \to \infty} \rho(x_{n_k}, z) = \tau(z).\end{aligned}$$

Therefore, from the uniqueness of $z$, we have $z = Tz$. □

**Remark 1** Theorem 4.8 in [14] is the special case of Theorem 1 when the space is a uniformly convex Banach space, that is, Theorem 1 generalizes Theorem 4.8 in [14]

**Lemma 2** Let $C$ be a nonempty closed convex subset of a uniformly convex partially ordered hyperbolic space $(E, \rho, \leq)$ and $T: C \to C$ be a monotone mapping on $C$. Assume that there exists $x_1 \in C$ such that $x_1 \leq Tx_1 (Tx_1 \leq x_1)$ and $\{x_n\}$ is a sequence generated by (**). Then

$$x_n \leq Tx_n \leq x_{n+1} (x_{n+1} \leq Tx_n \leq x_n)$$

holds for any $n \in \mathbf{N}$.

**Proof** Assume that . From the definition of (**) and the convexity of order intervals, we have

$$x_1 \leq z_1 \leq Tx_1 \leq Tz_1$$
$$\Rightarrow z_1 \leq y_1 \leq Tz_1 \leq Ty_1$$
$$\Rightarrow y_1 \leq Tz_1 \leq x_2 \leq Ty_1 \leq Tx_2.$$

Hence we have $x_1 \leq Tx_1 \leq x_2 \leq Tx_2$.

By inductive way, we obtain $x_n \leq Tx_n \leq x_{n+1}$ for all $n \in \mathbf{N}$. □

**Theorem 2** Let $C$ be a nonempty closed convex subset of a uniformly convex partially ordered hyperbolic space $(E, \rho, \leq)$ and $T: C \to C$ be a monotone generalized $\alpha$-noexpansive mapping on $C$. Assume that there exists $x_1 \in C$ such that $x_1 \leq Tx_1 (Tx_1 \leq x_1)$ and $\{x_n\}$ is a sequence generated by (**). Let $\{x_n\}$ be bounded and there exists $x \in C$ such that $x_n \leq x(x \leq x_n)$ for all $n \in \mathbf{N}$, and $\liminf_{n \to \infty} \rho(x_n, Tx_n) = 0$. Then $T$ has a fixed point.

**Proof.** From Lemma 2 and Theorem 1, we can prove the result easily. □

## 3. Convergence theorems

In this section, we establish $\Delta$ – convergence and strong convergence theorems of (**) for a monotone generalized $\alpha$-nonexpansive mapping.

First we give some useful lemmas.

**Lemma 3** ([7]) Let $E$ be a uniformly convex hyperbolic space with monotone modulus of uniform convexity and $0 < a \leq s_n \leq b < 1$ for all $n \geq 1$. And assume that $\{x_n\}$ and $\{y_n\}$ are sequences in $E$ such that

$$\exists d \geq 0, \limsup_{n \to \infty} \rho(x_n, z) \leq d, \limsup_{n \to \infty} \rho(y_n, z) \leq d,$$
$$\limsup_{n \to \infty} \rho(s_n x_n \oplus (1 - s_n) y_n, z) = d$$

for $z \in E$.

Then we have $\lim_{n \to \infty} \rho(x_n, y_n) = 0$.

**Lemma 4** Let $C$ be a nonempty closed convex subset of a partially ordered hyperbolic space $(E, \rho, \leq)$ and $T: C \to C$ be a monotone generalized $\alpha$- noexpansive mapping on $C$. Assume that there exists $x_1 \in C$ such that $x_1 \leq Tx_1 (Tx_1 \leq x_1)$ and $\{x_n\}$ is a sequence generated by (**). If $T$ has a fixed point in $C$ and $x_1 \leq p (p \leq x_1)$ for any fixed point $p$, then the limit $\lim_{n \to \infty} \rho(x_n, p)$ exists.

**Proof** Assume that $x_1 \leq Tx_1$. Since $T$ is monotone, we have

$$z_1 \leq Tx_1 \leq Tp = p$$
$$\Rightarrow y_1 \leq Tz_1 \leq Tp = p$$
$$\Rightarrow x_2 \leq Ty_1 \leq Tp = p$$
$$\Rightarrow x_2 \leq p$$

By inductive way, we have $x_n \leq p$ for all $n \in \mathbf{N}$. Since $T$ is monotone quasi-nonexpansive, we have

$$\begin{aligned} \rho(z_n, p) &= \rho((1-c_n)x_n \oplus c_n Tx_n, p) \\ &\leq (1-c_n)\rho(x_n, p) + c_n \rho(Tx_n, p) \\ &\leq (1-c_n)\rho(x_n, p) + c_n \rho(x_n, p) \\ &= \rho(x_n, p) \end{aligned} \quad (1)$$

and

$$\begin{aligned} \rho(y_n, p) &= \rho((1-b_n)z_n \oplus b_n Tz_n, p) \\ &\leq (1-b_n)\rho(z_n, p) + b_n \rho(Tz_n, p) \\ &\leq (1-b_n)\rho(z_n, p) + b_n \rho(z_n, p) = \rho(z_n, p) \\ &\leq \rho(x_n, p). \end{aligned} \quad (2)$$

From (1) and (2), we have

$$\begin{aligned} \rho(x_{n+1}, p) &= \rho((1-a_n)Tz_n \oplus a_n Ty_n, p) \\ &\leq (1-a_n)\rho(Tz_n, p) + a_n \rho(Ty_n, p) \\ &\leq (1-a_n)\rho(x_n, p) + a_n \rho(x_n, p) \\ &= \rho(x_n, p). \end{aligned}$$

Therefore, the sequence $\{\rho(x_n, p)\}$ is decreasing and bounded below. So the limit $\lim_{n \to \infty} \rho(x_n, p)$ exists. □

**Lemma 5** Let $C$ be a nonempty closed convex subset of a complete uniformly convex partially ordered hyperbolic space $(E, \rho, \leq)$ with monotone modulus of uniform convexity and $T: C \to C$ be a monotone generalized $\alpha$-noexpansive mapping on $C$. Assume that there exists $x_1 \in C$ such that $x_1 \leq Tx_1 (Tx_1 \leq x_1)$ and $\{x_n\}$ is a sequence generated by (**). If $T$ has a fixed point in $C$ and $x_1 \leq p (p \leq x_1)$ for any fixed point $p$, then $\lim_{n \to \infty} \rho(x_n, Tx_n) = 0$.

**Proof** Form lemma 4, $\lim_{n \to \infty} \rho(x_n, p)$ exists. Set

$$\lim_{n \to \infty} \rho(x_n, p) = d . \tag{3}$$

Then from (1), (2) and (3), we have

$$\limsup_{n \to \infty} \rho(z_n, p) \leq \limsup_{n \to \infty} \rho(x_n, p) \leq d$$

$$\limsup_{n \to \infty} \rho(y_n, p) \leq \limsup_{n \to \infty} \rho(x_n, p) \leq d .$$

Since a monotone generalized $\alpha$-noexpansive mapping is monotone quasi-nonexpansive, we have

$$\rho(Tx_n, p) \leq \rho(x_n, p) \Rightarrow \limsup_{n \to \infty} \rho(Tx_n, p) \leq \limsup_{n \to \infty} \rho(x_n, p) \leq d . \tag{4}$$

Similarly we have

$$\limsup_{n \to \infty} \rho(Ty_n, p) \leq \limsup_{n \to \infty} \rho(y_n, p) \leq d , \tag{5}$$

$$\limsup_{n \to \infty} \rho(Tz_n, p) \leq \limsup_{n \to \infty} \rho(z_n, p) \leq d . \tag{6}$$

Also we obtain

$$d = \lim_{n \to \infty} \rho(x_{n+1}, p) = \lim_{n \to \infty} \rho((1-a_n)Tz_n \oplus a_n Ty_n, p) \tag{7}$$

Then from (5)-(7) and lemma 3, we have

$$\lim_{n \to \infty} \rho(Tz_n, Ty_n) = 0 . \tag{8}$$

Since

$$\rho(x_{n+1}, Tz_n) = \rho((1-a_n)Tz_n \oplus a_n Ty_n, Tz_n)$$
$$\leq a_n \rho(Tz_n, Ty_n)$$

, from (8) we have $\lim_{n \to \infty} \rho(x_{n+1}, Tz_n) = 0$.

Also we have

$$\rho(x_{n+1}, p) \leq \rho(x_{n+1}, Tz_n) + \rho(Tz_n, p)$$
$$\leq \rho(x_{n+1}, Tz_n) + \rho(z_n, p)$$

and by taking limits on both sides, we have the following from (6).

$$d = \liminf_{n\to\infty} \rho(x_{n+1}, p) \leq \liminf_{n\to\infty} \rho(z_n, p)$$
$$\Rightarrow d = \lim_{n\to\infty} \rho(z_n, p)$$

Hence we obtain

$$d = \lim_{n\to\infty} \rho(z_n, p) = \lim_{n\to\infty} \rho((1-c_n)x_n \oplus c_n Tx_n, p) \tag{9}$$

Therefore, from (4), (9) and lemma 3, we have $\lim_{n\to\infty} \rho(x_n, Tx_n) = 0$. □

**Theorem 3** Let $C$ be a nonempty closed convex subset of a complete uniformly convex partially ordered hyperbolic space $(E, \rho, \leq)$ with monotone modulus of uniform convexity and $T: C \to C$ be a monotone generalized $\alpha$-noexpansive mapping on $C$. Assume that there exists $x_1 \in C$ such that $x_1 \leq Tx_1 (Tx_1 \leq x_1)$ and $\{x_n\}$ is a sequence generated by (**). If $T$ has a fixed point in $C$ and $x_1 \leq p (p \leq x_1)$ for any fixed point $p$, then $\{x_n\}$ $\Delta$-converges to a fixed point.

**Proof** From lemma 4, for any fixed point $p$, $\lim_{n\to\infty} \rho(x_n, p)$ exists. So $\{x_n\}$ is bounded and the asymptotic center of $\{x_n\}$ with respect to $C$ is a unique point $z_0 \in C$. Let $z$ be the asymptotic center of any subsequence $\{x_{n_k}\}$ of $\{x_n\}$ with respect to $C$. Now we prove that $z$ is a fixed point.

From lemma 5, $\lim_{n\to\infty} \rho(x_n, Tx_n) = 0$. Then from lemma 1, we have

$$\limsup_{k\to\infty} \rho(x_{n_k}, Tz) \leq \frac{3+\alpha}{1-\alpha} \limsup_{k\to\infty} \rho(x_{n_k}, Tx_n) + \limsup_{k\to\infty} \rho(x_{n_k}, z)$$
$$= \limsup_{k\to\infty} \rho(x_{n_k}, z).$$

By the uniqueness of an asymptotic center, we have $z = Tz$.

Now we prove that $z_0 = z$.

Assume that $z_0 \neq z$, then we have the following from the definition of an asymptotic center.

$$\limsup_{k\to\infty} \rho(x_{n_k}, z) < \limsup_{k\to\infty} \rho(x_{n_k}, z_0)$$
$$= \limsup_{n\to\infty} \rho(x_n, z_0)$$
$$< \limsup_{n\to\infty} \rho(x_n, z)$$
$$= \limsup_{k\to\infty} \rho(x_{n_k}, z)$$

This is a contradiction, so $z_0 = z$.

Therefore $\{x_n\}$ $\Delta$–converges to $z \in C$. $\square$

**Theorem 4** Let $C$ be a nonempty closed convex subset of a complete uniformly convex partially ordered hyperbolic space $(E, \rho, \leq)$ with monotone modulus of uniform convexity and $T: C \to C$ be a monotone generalized $\alpha$-noexpansive mapping on $C$. Assume that there exists $x_1 \in C$ such that $x_1 \leq Tx_1 (Tx_1 \leq x_1)$ and $\{x_n\}$ is a sequence generated by (**). If $T$ has a fixed point in $C$ and $x_1 \leq p (p \leq x_1)$ for any fixed point $p$, then $\{x_n\}$ converges to a fixed point if and only if

$$\liminf_{n\to\infty} d(x_n, F(T)) = 0$$

, where $d(x_n, F(T)) = \inf\{\rho(x_n, p) : p \in F(T)\}$.

**Proof** Necessity is clear.

Assume that $\liminf_{n\to\infty} d(x_n, F(T)) = 0$. First we prove that $F(T)$ is closed.

Assume that a sequence $\{z_n\}$ of $F(T)$ converges to $z \in C$.

Since $\frac{1}{2}\rho(z_n, Tz_n) = 0 < \rho(z_n, z)$ for all $n \in \mathbf{N}$, we have

$$\rho(Tz_n, Tz) = \rho(z_n, Tz) \leq \alpha\rho(z_n, Tz_n) + \alpha\rho(z, Tz) + (1-2\alpha)\rho(z_n, z)$$
$$\leq \alpha\rho(z_n, z) + \alpha\rho(z, Tz_n) + (1-2\alpha)\rho(z_n, z)$$
$$= \alpha\rho(Tz_n, z) + (1-\alpha)\rho(z_n, z).$$

That is $\rho(z_n, Tz) \leq \rho(z_n, z)$. So $z \in F(T)$ and $F(T)$ is closed.

Since $\liminf_{n\to\infty} d(x_n, F(T)) = 0$, for any sequence $\{z_k\}$ of $F(T)$, there exists a subsequence $\{x_{n_k}\}$ of $\{x_n\}$ such that for all $k \in \mathbf{N}$,

$$\rho(x_{n_k}, z_k) \leq \frac{1}{2^k}.$$

By triangle inequality, we have

$$\rho(z_{k+1}, z_k) \leq \rho(z_{k+1}, x_{n_{k+1}}) + \rho(x_{n_{k+1}}, z_k) \leq \frac{1}{2^{k+1}} + \frac{1}{2^k} < \frac{1}{2^{k-1}}.$$

Hence $\{z_k\}$ is Cauchy sequence, so it converges to $z \in F(T)$. By triangle inequality, we have

$$\rho(x_{n_k}, z) \leq \rho(x_{n_k}, z_k) + \rho(z, z_k).$$

By taking limits as $k \to \infty$, $\{x_{n_k}\}$ converges to $z \in F(T)$. From lemma 4 $\lim_{n \to \infty} \rho(x_n, p)$ exists, so we have

$$\lim_{n \to \infty} \rho(x_n, z) = \lim_{k \to \infty} \rho(x_{n_k}, z) = 0.$$

Therefore $\{x_n\}$ converges to $z \in F(T)$. □

Next we establish the strong convergence theorem of (**) by using ($I$)-condition.

**Theorem 5** Let $C$ be a nonempty closed convex subset of a complete uniformly convex partially ordered hyperbolic space $(E, \rho, \leq)$ with monotone modulus of uniform convexity and $T: C \to C$ be a monotone generalized $\alpha$-noexpansive mapping on $C$ satisfying the ($I$)-condition. Assume that there exists $x_1 \in C$ such that $x_1 \leq Tx_1 (Tx_1 \leq x_1)$ and $\{x_n\}$ is a sequence generated by (**). If $T$ has a fixed point in $C$ and $x_1 \leq p (p \leq x_1)$ for any fixed point $p$, then $\{x_n\}$ converges to a fixed point.

**Proof** From lemma 5, $\lim_{n \to \infty} \rho(x_n, Tx_n) = 0$. From the ($I$)-condition, we have

$$0 \leq \lim_{n \to \infty} h(d(x_n, F(T))) \leq \lim_{n \to \infty} \rho(x_n, Tx_n) = 0$$
$$\Rightarrow \lim_{n \to \infty} h(d(x_n, F(T))) = 0.$$

Since $h: [0, \infty) \to [0, \infty)$ decreases and $h(0) = 0, h(r) > 0, \forall r > 0$, we have

$$\lim_{n \to \infty} d(x_n, F(T)) = 0.$$

All the assumptions in theorem 4 are satisfied, so $\{x_n\}$ converges to a fixed point of $T$. □

**Remark 2** Theorem 5 implies the existence of the minimum(maximum) fixed point and the sequence converging to this.

## 4. Numerical example

In this section, we show that (**) converges to a fixed point faster than (**) by a numerical example.

Let the mapping $T:[0, 4] \to [0, 4]$ be defined as follows.

$$T(x) = \begin{cases} 0, & x \neq 4, \\ 2, & x = 4. \end{cases}$$

$T$ is not a ($C$)-mapping but a monotone generalized $\alpha$-nonexpansive mapping for $\alpha \geq \frac{1}{3}$ ([14]). And it has a unique fixed point 0.

With help of Mathematica Program Software, we obtain the comparison Table 1 and Figure 1 for (*) and (**) with $a_n = 0.85, b_n = 0.65, c_n = 0.45, x_1 = 0.9$.

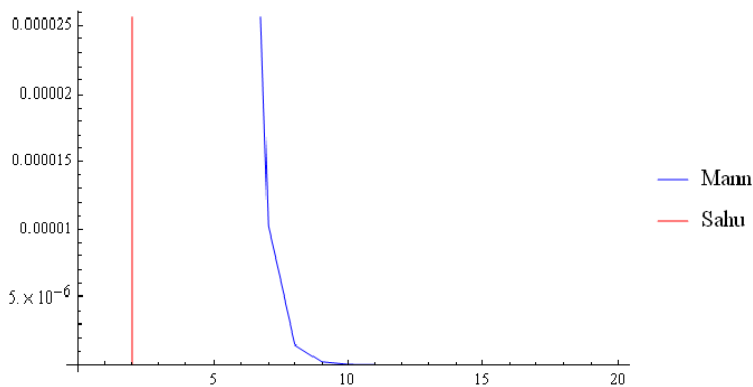

Figure 1

|    | Mann                     | Sahu |
|----|--------------------------|------|
| 1  | 0.9                      | 0.9  |
| 2  | 0.135                    | 0.   |
| 3  | 0.02025                  | 0.   |
| 4  | 0.0030375                | 0.   |
| 5  | 0.000455625              | 0.   |
| 6  | 0.0000683438             | 0.   |
| 7  | 0.0000102516             | 0.   |
| 8  | $1.53773 \times 10^{-6}$ | 0.   |
| 9  | $2.3066 \times 10^{-7}$  | 0.   |
| 10 | $3.4599 \times 10^{-8}$  | 0.   |
| 11 | $5.18985 \times 10^{-9}$ | 0.   |
| 12 | $7.78478 \times 10^{-10}$| 0.   |
| 13 | $1.16772 \times 10^{-10}$| 0.   |
| 14 | $1.75158 \times 10^{-11}$| 0.   |
| 15 | $2.62736 \times 10^{-12}$| 0.   |
| 16 | $3.94105 \times 10^{-13}$| 0.   |
| 17 | $5.91157 \times 10^{-14}$| 0.   |
| 18 | $8.86735 \times 10^{-15}$| 0.   |
| 19 | $1.3301 \times 10^{-15}$ | 0.   |
| 20 | $1.99515 \times 10^{-16}$| 0.   |

Table 1

From the figure and table, we can see that (**) converges faster than (*) to the fixed point $p = 0$.

## 5. Application to integral equations

In this section, we apply the results in previous sections to solve the following inequality equation in Hilbert space $L^2([0, 1], \mathbf{R})$.

$$x(t) = y_0(t) + \int_0^1 B(t, z, x(z))dz, t \in [0, 1] \qquad (10)$$

Here are some basic concepts.

**Definition 9** ([3]) Let $P$ be a subset of a Banach space $E$. $P$ is said to be an order cone if it satisfies the following conditions.

(i) $P$ is nonempty and $P \neq \{\theta\}$;

(ii) $ax + by \in P$, for any $x, y \in P$ and any nonnegative real numbers $a, b$;

(iii) $P \cap (-P) = \{\theta\}$.

If $P \subseteq E$ is an order cone, then the partial order in Banach space $E$ with respect to $P$ can be defined as follows.

$$x \leq y \overset{def}{\Longleftrightarrow} y - x \in P.$$

**Definition 10** ([3]) Let $P$ be the order cone of the partially ordered Banach space $E$. $P$ is said to be fully regular if every increasing sequence of $E$ which is norm bounded converges.

That is, for the sequence $\{a_n\}$ of $E$ satisfying

$$a_1 \leq a_2 \leq \cdots \leq a_n \leq \cdots,$$
$$\exists M, \forall n \in \mathbf{N} : \|a_n\| \leq M$$

, there exists $a \in E$ such that $\lim_{n \to \infty} \|a_n - a\| = 0$.

The order cone of $L^p([0,1], \mathbf{R})(1 \leq p < \infty)$ is fully regular.

With sufficiently large $r$, set a closed convex subset $C$ as follows.

$$C := \{x \in L^2([0,1], \mathbf{R}) : \|x\| \leq r\}.$$

And define the mapping $T$ on $L^2([0,1], \mathbf{R})$ as follows.

$$Tx(t) := y_0(t) + \int_0^1 B(t, z, x(z)) dz.$$

The fully regular cone $P = \{g \in L^2([0,1], \mathbf{R}) : g(t) \geq 0, a.e.\, t \in [0,1]\}$ define the following partial order in $L^2([0,1], \mathbf{R})$.

$$u, v \in L^2([0,1], \mathbf{R}), u \leq v := u(t) \leq v(t), a.e.\, t \in [0,1]$$

**Lemma 6** Assume that the following assumptions hold.

(i) $y_0 \in L^2([0,1], \mathbf{R})$;

(ii) $B : [0,1] \times [0,1] \times L^2([0,1], \mathbf{R}) \to \mathbf{R}$ is nonnegative, continuous and satisfies the following two properties for any $t, z \in [0,1]$.

(a) For any $u, v \in L^2([0,1], \mathbf{R})$ with $u \leq v$,

$$0 \leq B(t, z, v) - B(t, z, u) \leq \|v - u\|. \tag{11}$$

(b) There exists a nonnegative function $f(\cdot, \cdot) \in L^2([0,1] \times [0,1])$ and $M < \dfrac{1}{2}$ such that for any $x \in L^2([0,1], \mathbf{R})$,

$$|B(t, z, x)| \leq f(t, z) + M\|x\| \tag{12}$$

Then $T$ is monotone nonexpansive such that $T(C) \subset C$ and $y_0 \leq Tx$ for any $x \in C$.

**Proof** From the nonnegativeness of $B$, we have $y_0 \leq Tx$ for any $x \in C$.

Now we prove that $T(C) \subset C$.

For $x \in C$, from (12) we have

$$\|Tx(t)\| = \left\| y_0(t) + \int_0^1 B(t, z, x(z))dz \right\|$$

$$\leq \|y_0\| + \left\| \int_0^1 B(t, z, x(z))dz \right\|$$

$$\leq \|y_0\| + \left\| \int_0^1 f(t, z)dz \right\| + M\|x\|$$

So by taking $r$ such that $\frac{1}{2}r \geq \|y_0\| + \left\| \int_0^1 f(t, z)dz \right\|$, we have $Tx \in C$.

Next we prove that $T$ is monotone.

From (11), for any $u, v \in C$ with $u \leq v$, we have

$$Tu = y_0(t) + \int_0^1 B(t, z, u(z))dz \leq y_0(t) + \int_0^1 B(t, z, v(z))dz = Tv.$$

Therefore $T$ is monotone.

Finally we prove the nonexpansiveness.

From (11), we have

$$\|Tv - Tu\| = \left\| \int_0^1 (B(t, z, v) - B(t, z, u))dz \right\| \leq \int_0^1 \|v - u\|dz = \|v - u\|$$

, for any $u, v \in C$ with $u \leq v$.

This completes the proof. $\square$

Now we prove the existence of a solution of the equation (10).

**Theorem 5** Under the assumptions in lemma 6, $T$ has a fixed point. That is, (10) has a solution.

**Proof** Set $x_1 := y_0$. From lemma 6, $x_1 \leq Tx_1$.

For all $n \in \mathbf{N}$, set $x_{n+1} := Tx_n$. Since $T$ is monotone, we have $x_n \leq x_{n+1}$ for all $n \in \mathbf{N}$. Hence $\{x_n\}$ is increasing.

Since the order cone is fully regular, $\{x_n\}$ converges to its supremum $x$. From the nonexpansiveness, $x$ is a fixed point of $T$. This completes the proof. □

Now we establish the convergence theorem for the solution of (10) by using the convergence result in section 3.

We can easily prove that a Hilbert space is a uniformly convex hyperbolic metric space with monotone modulus of uniform convexity. In fact, since

$$\|x+y\|^2 = 2(\|x\|^2 + \|y\|^2) - \|x-y\|^2$$

for any $x, y$ in a Hilbert space, we have

$$\|x+y\|^2 \leq (4-\varepsilon^2)r^2$$

$$\Rightarrow \left\|\frac{x+y}{2}\right\| \leq r \cdot \sqrt{1-\frac{\varepsilon^2}{4}}$$

if $\|x\| \leq r, \|y\| \leq r, \|x-y\| \geq r\varepsilon$.

Therefore the modulus of uniform convexity is

$$\delta(r,\varepsilon) = 1 - \sqrt{1-\frac{\varepsilon^2}{4}}.$$

Hence the modulus of uniform convexity is monotone.

We obtain the following result from theorem 2 in section 3.

**Theorem 6** Under the assumptions in lemma 6, let $\{x_n\}$ be a sequence by (\*\*) with $x_1 = y_0$. Then $\{x_n\}$ $\Delta$-converges to a solution of (10).